\documentclass[11pt]{article}
\usepackage{latexsym}
\newcommand{\eproof}{\mbox{\ }\hfill $\Box$ \par \vskip 10pt}

\newtheorem{Theorem}{Theorem}[section]
\newtheorem{lemma}[Theorem]{Lemma}
\newtheorem{prop}[Theorem]{Proposition}

\newtheorem{corol}[Theorem]{Corollary}

\begin{document}

\title{Boundary stabilization of transmission problems}

\author{{\sc Fernando Cardoso and Georgi Vodev\thanks{corresponding author}}}

\date{}

\maketitle

\noindent
{\bf Abstract.} We study the transmission problem in bounded domains with dissipative boundary conditions. Under some natural assumptions, we prove  uniform bounds of the corresponding resolvents on the real axis at high frequency, and as a consequence, we obtain free of eigenvalues regions. 
To this end, we extend the result of \cite{kn:CPV} under more general assumptions. 
As an application, we get exponential decay of the energy of the solutions of the corresponding mixed boundary value problems.

\setcounter{section}{0}
\section{Introduction and statement of results}

Let $\Omega_1\subset\Omega_2\subset...\subset\Omega_{m+1}\subset {\bf R}^n$, $m\ge 1$, $n\ge 2$, be bounded, strictly convex domains with smooth boundaries $\Gamma_k=\partial\Omega_k$, $\Gamma_k\cap\Gamma_{k+1}=\emptyset$. Let also $\Omega_0\subset\Omega_1$ be a bounded domain with smooth boundary $\Gamma_0=\partial\Omega_0$ such that ${\bf R}^n\setminus\Omega_0$ is connected. In the present paper we are interested in studying the large time behavior of the solutions of the following mixed boundary value problems:
$$
\left\{
\begin{array}{l}
 (i\partial_t-c_k^2\Delta)u_k^0(x,t)=0\quad\mbox{in}\quad(\Omega_k\setminus\Omega_{k-1})\times (0,+\infty),\,k=1,...,m+1,\\
Bu_1^0(x,t)=0\quad\mbox{on}\quad\Gamma_0\times (0,+\infty),\\
u_k^0(x,t)=u_{k+1}^0(x,t),\,\partial_\nu u_k^0(x,t)=\partial_\nu u_{k+1}^0(x,t)
\quad\mbox{on}\quad \Gamma_k\times(0,+\infty),\,k=1,...,m,\\
\partial_\nu u_{m+1}^0(x,t)+ia(x)u_{m+1}^0(x,t)=0\quad\mbox{on}\quad \Gamma_{m+1}\times(0,+\infty),\\
u_k^0(x,0) = f_k^0(x),  \, k=1,...,m+1,
\end{array}
\right.
\eqno{(1.1)}
$$
and
$$
\left\{
\begin{array}{l}
 (\partial_t^2-c_k^2\Delta)u_k^1(x,t)=0\quad\mbox{in}\quad(\Omega_k\setminus\Omega_{k-1})\times (0,+\infty),\,k=1,...,m+1,\\
Bu_1^1(x,t)=0\quad\mbox{on}\quad\Gamma_0\times (0,+\infty),\\
u_k^1(x,t)=u_{k+1}^1(x,t),\,\partial_\nu u_k^1(x,t)=\partial_\nu u_{k+1}^1(x,t)
\quad\mbox{on}\quad \Gamma_k\times(0,+\infty),\,k=1,...,m,\\
\partial_\nu u_{m+1}^1(x,t)+a(x)\partial_tu_{m+1}^1(x,t)=0\quad\mbox{on}\quad \Gamma_{m+1}\times(0,+\infty),\\
u_k^1(x,0) = f_k^1(x), \,  \partial_t u_k^1(x,0) = g_k^1(x), \, k=1,...,m+1,
\end{array}
\right.
\eqno{(1.2)}
$$
where either $B=Id$ (Dirichlet boundary conditions) or $B=\partial_\nu$ (Neumann boundary conditions), $\partial_\nu$ denotes the normal derivative to the boundary, $c_k$ are constants satisfying
$$c_1>c_2>...>c_{m+1}>0,\eqno{(1.3)}$$
and $a(x)$ is a continuous, real-valued function on $\Gamma_{m+1}$ supposed to satisfy
$$a(x)\ge a_0 \quad\mbox{on}\quad \Gamma_{m+1},\eqno{(1.4)}$$
with some constant $a_0>0$. The equation (1.2) describes the propagation of acoustic waves in different media with different speeds $c_k$, $k=1,...,m+1$, which do not penetrate into $\Omega_0$.
The boundary condition on $\Gamma_{m+1}$ is a strong dissipative one which guarantees that the energy of the solutions of (1.2) with finite energy initial data tends to zero as $t\to +\infty$. The equation (1.1) is of Schr\"odinger type with weak dissipative boundary conditions. In fact, the large time behavior of the solutions of (1.1) and (1.2) is closely related to the behavior on the real axis of the corresponding resolvent operator, $R^j(\lambda)$, $\lambda\in{\bf C}$, defined for ${\rm Im}\,\lambda<0$ as follows. Given $$v^j=(v_1^j,...,v_{m+1}^j)\in H:=\oplus_{k=1}^{m+1}L^2\left(\Omega_{k}\setminus\Omega_{k-1},c_k^{-2}dx\right),$$
$R^j(\lambda)v^j=(u_1^j,...,u_{m+1}^j)\in H$ solves the equation
$$
\left\{
\begin{array}{l}
 (\lambda^2+c_k^2\Delta)u_k^j=v_k^j\quad\mbox{in}\quad\Omega_k\setminus\Omega_{k-1},\,k=1,...,m+1,\\
Bu_1^j=0\quad\mbox{on}\quad\Gamma_0,\\
u_k^j=u_{k+1}^j,\,\partial_\nu u_k^j=\partial_\nu u_{k+1}^j
\quad\mbox{on}\quad \Gamma_k,\,k=1,...,m,\\
\partial_\nu u_{m+1}^j+i\lambda^j a(x) u_{m+1}^j=0\quad\mbox{on}\quad \Gamma_{m+1}.
\end{array}
\right.
\eqno{(1.5)}
$$
It is well known that $\lambda^jR^j(\lambda):H\to H$ extends meromorphically to the whole complex plane ${\bf C}$ with no poles on the real axis (the latter can be derived from the Carleman estimates of \cite{kn:B}). In the present paper we will study the behavior of $R^j(\lambda)$
for $\lambda\in{\bf R}$, $|\lambda|\gg 1$. To this end we need to impose some conditions on $\Omega_0$ (as weak as possible). We first make the following assumption:
$$\mbox{every generalized ray in $\Omega_1\setminus\Omega_0$ hits the boundary $\Gamma_1$}.\eqno{(1.6)}$$
Clearly, (1.6) is fulfilled if $\Omega_0$ is strictly convex. However, the class of the domains for which (1.6) is satisfied is much larger than the class of strictly convex domains. We can now state our first result.

\begin{Theorem} Assume (1.3), (1.4) and (1.6) fulfilled. Then, there exist constants $C,C_1>0$ so that $R^j(\lambda)$ ($j=0,1$) satisfies the bound
$$\left\|R^j(\lambda)\right\|_{H\to H}\le C|\lambda|^{-j},\quad\lambda\in{\bf R},\,|\lambda|\ge C_1.\eqno{(1.7)}$$
\end{Theorem}

One can derive from this theorem the following

\begin{corol} Under the assumptions of Theorem 1.1, the solutions $$u^j(x,t)=(u_1^j(x,t),...,u_{m+1}^j(x,t))$$ of (1.1) and (1.2) satisfy the estimates (for $t\gg 1$):
$$\left\| u^0(\cdot,t)\right\|_H\le \widetilde Ce^{-Ct}\left\|u^0(\cdot,0)\right\|_H,\eqno{(1.8)}$$
with constants $\widetilde C,C>0$ independent of $t$ and $u^0$, and
$$\left\|\nabla_x u^1(\cdot,t)\right\|_H+\left\|\partial_t u^1(\cdot,t)\right\|_H\le \widetilde Ce^{-Ct}\left(\left\|\nabla_x u^1(\cdot,0)\right\|_H+\left\|\partial_t u^1(\cdot,0)\right\|_H\right),\eqno{(1.9)}$$
with constants $\widetilde C,C>0$ independent of $t$ and $u^1$.  
\end{corol}

To prove (1.8) and (1.9) it suffices to show that the solutions of (1.1) and (1.2) are given by semi-groups $e^{itA_j}$, respectively, acting on  suitable Hilbert spaces 
${\cal H}_j$ with generators $A_j$ of compact resolvent and hence of discrete spectrum. Then Theorem 1.1 implies that 
$$\left\|(A_j-z)^{-1}\right\|_{{\cal H}_j\to{\cal H}_j}=O(1)\quad\mbox{for}\quad z\in{\bf R},\,|z|\gg 1,$$
which in turn implies (1.8) and (1.9), respectively (see Section 2 for more details).

In the case when there is no transmission of waves (which corresponds to taking $m=0$ in the setting above) the above estimates follow from the results of \cite{kn:BLR}. In fact, in \cite{kn:BLR} a more general situation is studied, namelly $\Omega_1$ is not necessarilly strictly convex and (1.4) is supposed to hold on a non-empty subset $\widetilde\Gamma_1$ of $\Gamma_1$. Then (1.6) is replaced by the assumption that every generalized ray in $\Omega_1\setminus\Omega_0$ hits $\widetilde\Gamma_1$ at a non-diffractive point (see \cite{kn:BLR} for the definition and more details). The situation changes drastically in the case of transmission (which corresponds to taking $m\ge 1$ in the setting above) due to the fact that the classical flow for this problem is much more complicated. Indeed, when a ray in $\Omega_{k+1}\setminus\Omega_k$ hits the boundary $\Gamma_k$ (if $1\le k\le m$) or the boundary $\Gamma_{k+1}$ (if $0\le k\le m-1$), it splits into two rays - one staying in 
$\Omega_{k+1}\setminus\Omega_k$ and another entering into $\Omega_{k}\setminus\Omega_{k-1}$ or $\Omega_{k+2}\setminus\Omega_{k+1}$, respectively. Consequently, there are infinitely many rays which do not reach the boundary $\Gamma_{m+1}$ where the dissipation is active. The condition (1.3), however, guarantees that these rays carry a negligible amount of energy, and therefore (1.3) is crucial for the above estimates to hold. Indeed, if for example we have $c_{k_0}<c_{k_0+1}$
for some $1\le k_0\le m$, then one can construct quasi-modes concentrated on the boundary $\Gamma_{k_0}$ (see \cite{kn:PV}). Consequently, we have in this case a sequence, $\{\lambda_k\}_{k=1}^\infty$, of poles of $R^j(\lambda)$ such that $|\lambda_k|\to \infty$ and $0<{\rm Im}\,\lambda_k\le C_N|\lambda_k|^{-N}$, $\forall N\ge 1$. Note also that the fact that the domains $\Omega_k$, $k=1,...,m+1$, are strictly convex is crucial for our proof to work, and
quite probably Theorem 1.1 as well as the estimates (1.8) and (1.9) are no longer true without this condition. This is essential for the proof of Proposition 2.3 below (proved in \cite{kn:CPV}). It also guarantees nice properties of the Neumann operator (denoted by $N_k(\lambda)$, $k=1,...,m$ below) associated to the Helmholtz equation in ${\bf R}^n\setminus \Omega_k$ (see Lemmas 4.2 and 4.4).

To prove Theorem 1.1 we make use of the results of \cite{kn:CPV} where an exterior transmission problem has been studied. Consider the exterior stationary problem
$$
\left\{
\begin{array}{l}
 (\lambda^2+c_1^2\Delta)u=v\quad\mbox{in}\quad{\bf R}^n\setminus\Omega_0,\\
Bu=0\quad\mbox{on}\quad\Gamma_0,\\
u - \lambda - \mbox{outgoing}.
\end{array}
\right.
\eqno{(1.10)}
$$
Then the outgoing resolvent, ${\cal R}_0(\lambda)$, for the exterior problem is defined by $u={\cal R}_0(\lambda)v$. Let $\chi\in
C_0^\infty({\bf R}^n)$, $\chi=1$ on $\Omega_0$. It is well known that the cut-off resolvent $\chi {\cal R}_0(\lambda)\chi$ is analytic in ${\rm Im}\,\lambda<0$ and meromorphic in ${\rm Im}\,\lambda>0$ with no poles on the real axis. Clearly, the condition (1.6) implies that $\Omega_0$
is non-trapping, that is, all generalized rays in ${\bf R}^n\setminus\Omega_0$ escape at infinity. In particular, this implies that the cut-off resolvent $\chi {\cal R}_0(\lambda)\chi$ satisfies the bound
$$\left\|\chi {\cal R}_0(\lambda)\chi\right\|_{L^2({\bf R}^n\setminus\Omega_0)\to L^2({\bf R}^n\setminus\Omega_0)}\le C|\lambda|^{-1}\quad\mbox{for}\quad \lambda\in{\bf R},\,|\lambda|\ge 1.\eqno{(1.11)}$$
In fact, the only thing we use in the proof of Theorem 1.1 is the estimate (1.11). In other words, we can replace the condition (1.6) by the estimate (1.11). Note also that (1.11) implies that $\chi {\cal R}_0(\lambda)\chi$ extends analytically in a strip $\{\lambda\in{\bf C}:|{\rm Im}\,\lambda|\le Const,\,|\lambda|\ge 1\}$ and that (1.11) still holds in this larger region (see \cite{kn:V}).

An interesting open problem is to get estimates similar to those stated above for more general domains $\Omega_0$ for which (1.6) and (1.11)
are not satisfied. A typical example for such domains is $\Omega_0={\cal O}_1\cup{\cal O}_2$, where ${\cal O}_1$ and ${\cal O}_2$ are strictly convex domains with smooth boundaries, ${\cal O}_1\cap{\cal O}_2=\emptyset$. In this case there is one periodic ray between ${\cal O}_1$ and ${\cal O}_2$ which does not reach $\Gamma_1$. It is well known that in this case (1.11) does not hold. Instead, we have that, in the case of Dirichlet boundary conditions (i.e. $B=Id$), the cut-off resolvent $\chi {\cal R}_0(\lambda)\chi$ is analytic in a strip $\{\lambda\in{\bf C}:|{\rm Im}\,\lambda|\le Const,\,|\lambda|\ge 1\}$ with polynomially bounded norm (see \cite{kn:G}, \cite{kn:I1}). Our purpose is to treat such more general domains $\Omega_0$. More precisely, we make the following assumption:\\

There exist constants $C,C_1, C_2,p>0$ so that the cutoff resolvent $\chi {\cal R}_0(\lambda)\chi$ is analytic in a strip $\{\lambda\in{\bf C}:|{\rm Im}\,\lambda|\le C_1,\,|\lambda|\ge C_2\}$ and satisfies there the bound
 $$\left\|\chi {\cal R}_0(\lambda)\chi\right\|_{L^2({\bf R}^n\setminus\Omega_0)\to L^2({\bf R}^n\setminus\Omega_0)}\le C|\lambda|^p.\eqno{(1.12)}$$
 
 Note that (1.12) is also satisfied for domains $\Omega_0=\cup_{\ell=1}^L{\cal O}_\ell$, $L\ge 3$, where ${\cal O}_\ell$ are strictly convex domains with smooth boundaries, ${\cal O}_{\ell_1}\cap{\cal O}_{\ell_2}=\emptyset$, $\ell_1\neq\ell_2$, satisfying some natural conditions (see \cite{kn:I2} for more details). Note that in this case there could be infinitely many periodic broken rays which do not reach the boundary $\Gamma_1$. Let us also mention that semi-classical analogues of (1.12) have been recently proved in \cite{kn:NZ}, \cite{kn:NZ1} in a very general situation.
 
 Our main result is the following

\begin{Theorem} Assume (1.3), (1.4) and (1.12) fulfilled. Then, there exist constants $C,C_1>0$ so that $R^j(\lambda)$ ($j=0,1$) satisfies the bound
$$\left\|R^j(\lambda)\right\|_{H\to H}\le C|\lambda|^{-j}(\log|\lambda|)^{2^{m+1}},\quad\lambda\in{\bf R},\,|\lambda|\ge C_1.\eqno{(1.13)}$$ 
\end{Theorem}

Given an integer $k\ge 0$, set $\alpha_k=(2^k+1)^{-1}$. One can derive from this theorem the following

\begin{corol} Under the assumptions of Theorem 1.3, the solutions $$u^j(x,t)=(u_1^j(x,t),...,u_{m+1}^j(x,t))$$ of (1.1) and (1.2) satisfy the estimates (for $t\gg 1$):
$$\left\| u^0(\cdot,t)\right\|_H\le \widetilde C\exp\left(-C\varepsilon t^{\alpha_{m+1}}\right)\left\|u^0(\cdot,0)\right\|_{H^\varepsilon},\eqno{(1.14)}$$
for every $0<\varepsilon\le\varepsilon_0$, with constants $C,\varepsilon_0>0$ independent of $t$, $\varepsilon$ and $u^0$, $\widetilde C$ independent of $t$ and $u^0$, and
$$\left\|\nabla_x u^1(\cdot,t)\right\|_H+\left\|\partial_t u^1(\cdot,t)\right\|_H\le \widetilde C\exp\left(-C\varepsilon t^{\alpha_{m+1}}\right)\left(\left\|\nabla_x u^1(\cdot,0)\right\|_{H^\varepsilon}+\left\|\partial_t u^1(\cdot,0)\right\|_{H^\varepsilon}\right),\eqno{(1.15)}$$
for every $0<\varepsilon\le\varepsilon_0$, with constants $C,\varepsilon_0>0$ independent of $t$, $\varepsilon$ and $u^1$, $\widetilde C$ independent of $t$, and $u^1$,
where $H^\varepsilon :=\oplus_{k=1}^{m+1}H^\varepsilon\left(\Omega_{k}\setminus\Omega_{k-1}\right)$ denotes the corresponding Sobolev space.
\end{corol}

Note that the estimate (1.15) (with $\alpha_{m+1}=1/2$) has been proved in \cite{kn:C}, \cite{kn:C1} in the case of the damped wave equation on a bounded manifold without boundary under the assumption that there is only one closed geodesic of hyperbolic type which does not pass through the support of the dissipative term but 
all other geodesics do so. This result has been recently improved in \cite{kn:S} for a class of manifolds with negative curvature, where a strip free of eigenvalues has been obtained and, as a consequence, an analogue of (1.15) (with $\alpha_{m+1}=1$) has been proved. 

If $\Omega_0$ is strictly convex, the conclusions of Theorem 1.1 still hold if we admit transmision of waves in the interior of $\Omega_0$ moving with a speed $>c_1$, i.e. if we replace the boundary condition $Bu=0$ on $\Gamma_0$ by a transmission problem. Indeed, in this case we have (1.11) according to the results of \cite{kn:CPV}. Thus, it is natural to ask whether Theorem 1.3 still holds if $\Omega_0$ consists of two strictly convex bodies and we admit transmision of waves in the interior. To be more precise, we define the resolvent $\widetilde {\cal R}_0(\lambda)$ as $u=\widetilde{\cal R}_0(\lambda)v$, where $u=(u_1,u_2,u_3)$ and $v=(v_1,v_2,v_3)$ satisfy the equation
$$
\left\{
\begin{array}{l}
(\lambda^2+\alpha_k^2\Delta)u_k=v_k\quad\mbox{in}\quad{\cal O}_k,\,k=1,2,\\
 (\lambda^2+c_1^2\Delta)u_3=v_3\quad\mbox{in}\quad{\bf R}^n\setminus({\cal O}_1\cup{\cal O}_2),\\
u_k=u_3,\,\partial_\nu u_k=\partial_\nu u_3\quad\mbox{on}\quad\partial{\cal O}_k,\,k=1,2,\\
u_3 - \lambda - \mbox{outgoing},
\end{array}
\right.
\eqno{(1.16)}
$$
where $\alpha_k>c_1$, $k=1,2$, are constants, ${\cal O}_1$ and ${\cal O}_2$ are strictly convex domains with smooth boundaries, ${\cal O}_1\cap{\cal O}_2=\emptyset$. In analogy with the case of one strictly convex body discussed above, it is natural to make the following\\

\noindent
{\bf Conjecture.} {\it The resolvent} $\widetilde {\cal R}_0(\lambda)$ {\it satisfies the condition (1.12).}\\

Clearly, if this conjecture holds true, so does Theorem 1.3 in this more complex situation. However, it seems quite hard to prove.

The method we develop to prove the above results allows to get a decay of the local energy of the solutions of the following problem:
$$
\left\{
\begin{array}{l}
 (\partial_t^2-c_k^2\Delta)u_k(x,t)=0\quad\mbox{in}\quad(\Omega_k\setminus\Omega_{k-1})\times (0,+\infty),\,k=1,...,m,\\
 (\partial_t^2-c_{m+1}^2\Delta)u_{m+1}(x,t)=0\quad\mbox{in}\quad({\bf R}^n\setminus\Omega_m)\times (0,+\infty),\\
Bu_1(x,t)=0\quad\mbox{on}\quad\Gamma_0\times (0,+\infty),\\
u_k(x,t)=u_{k+1}(x,t),\,\partial_\nu u_k(x,t)=\partial_\nu u_{k+1}(x,t)
\quad\mbox{on}\quad \Gamma_k\times(0,+\infty),\,k=1,...,m,\\
u_k(x,0) = f_k(x), \,  \partial_t u_k(x,0) = g_k(x), \, k=1,...,m+1.
\end{array}
\right.
\eqno{(1.17)}
$$
More precisely, we have the following

\begin{Theorem} Under the assumptions (1.3) and (1.6), for every compact $K\subset{\bf R}^n\setminus\Omega_0$ there exists a constant 
$C_K>0$ so that the solution $$u(x,t)=(u_1(x,t),...,u_{m+1}(x,t))$$ of (1.17) satisfies the estimate (for $t\gg 1$)
$$\left\|\nabla_x u(\cdot,t)\right\|_{L^2(K)}+\left\|\partial_t u(\cdot,t)\right\|_{L^2(K)}\le C_Kp_0(t)\left(\left\|\nabla_x u(\cdot,0)\right\|_{L^2(K)}+\left\|\partial_t u(\cdot,0)\right\|_{L^2(K)}\right),\eqno{(1.18)}$$
provided ${\rm supp}\,u(\cdot,0)$, ${\rm supp}\,\partial_tu(\cdot,0)\subset K$, where
$$p_0(t)=\left\{\begin{array}{l} 
e^{-\gamma t}\quad if\,\,n\,\, is\,\, odd,\\
t^{-n}\quad if\,\,n\,\,is\,\, even,
\end{array}\right.$$
with a constant $\gamma>0$ independent of $t$. Furthermore, under the assumptions (1.3) and (1.12), we have the weaker estimate
$$\left\|\nabla_x u(\cdot,t)\right\|_{L^2(K)}+\left\|\partial_t u(\cdot,t)\right\|_{L^2(K)}\le C_{K,\varepsilon}p_\varepsilon(t)\left(\left\|\nabla_x u(\cdot,0)\right\|_{H^\varepsilon(K)}+\left\|\partial_t u(\cdot,0)\right\|_{H^\varepsilon(K)}\right),\eqno{(1.19)}$$
for every $0<\varepsilon\le\varepsilon_0$, provided ${\rm supp}\,u(\cdot,0)$, ${\rm supp}\,\partial_tu(\cdot,0)\subset K$, where
$$p_\varepsilon(t)=\left\{\begin{array}{l} 
\exp\left(-\gamma\varepsilon t^{\alpha_{m}}\right)\quad if\,\,n\,\,is\,\,odd,\\
t^{-n}\quad if\,\,n\,\, is\,\, even,
\end{array}\right.$$
with constants $\varepsilon_0,\gamma>0$ independent of $t$ and $\varepsilon$.
\end{Theorem}

Note that the estimate (1.18) is known to hold for non-trapping compactly supported perturbations of the Euclidean Laplacian (see \cite{kn:Va}). Note also that an estimate similar to (1.19) (with $\alpha_m=1/2$) has been proved in \cite{kn:C2} in the case of compactly supported metric perturbations of the Euclidean Laplacian under the assumption that there is only one closed geodesics of hyperbolic type. 

According to the results of \cite{kn:V}, to prove (1.18) it suffices to show that the corresponding cutoff resolvent is analytic in some
strip near the real axis with a suitable control of its norm at high frequencies. Thus, (1.18) follows from Theorem 2.2 below applied with $k=m$ (which is actually proved in \cite{kn:CPV}), while (1.19) is a consequence of Theorem 3.2 applied with $k=m$.

The paper is organized as follows. In Section 2 we prove Theorem 1.1, Corollary 1.2 and (1.18) using in an essential way the results of \cite{kn:CPV}. Similar ideas have already been used in \cite{kn:AV}. In Section 3 we prove Theorem 1.3, Corollary 1.4 and (1.19). To this end
we prove in Section 4 an analogue of the results of \cite{kn:CPV} under (1.12) (see Theorem 3.2 below).

\section{The case $\Omega_0$ non-trapping}

Let $w=(w_1,...,w_{m+1})$, $v=(v_1,...,v_{m+1})$ satisfy the equation
$$
\left\{
\begin{array}{l}
 (\lambda^2+c_k^2\Delta)w_k=v_k\quad\mbox{in}\quad\Omega_k\setminus\Omega_{k-1},\,k=1,...,m+1,\\
Bw_1=0\quad\mbox{on}\quad\Gamma_0,\\
w_k=w_{k+1},\,\partial_\nu w_k=\partial_\nu w_{k+1}
\quad\mbox{on}\quad \Gamma_k,\,k=1,...,m.
\end{array}
\right.
\eqno{(2.1)}
$$
We will first show that Theorem 1.1 follows from the following

\begin{Theorem} Assume (1.3) and (1.6) fulfilled. Then, there exist constants $C,\lambda_0>0$ so that for $\lambda\ge\lambda_0$ the solution to (2.1) satisfies the estimate
$$\|w\|_H\le C\lambda^{-1}\|v\|_H+C\left\|w_{m+1}|_{\Gamma_{m+1}}\right\|_{L^2(\Gamma_{m+1})}+C\lambda^{-1}\left\|\partial_\nu w_{m+1}|_{\Gamma_{m+1}}\right\|_{L^2(\Gamma_{m+1})}.\eqno{(2.2)}$$
\end{Theorem}

Applying Green's formula to the solution of (1.5) in each domain $\Omega_k\setminus\Omega_{k-1},\,k=1,...,m+1,$ and summing up these identities lead to the identity
$${\rm Im}\,\left\langle u^j,v^j\right\rangle_H:=\sum_{k=1}^{m+1}{\rm Im}\,\left\langle c_k^{-2}u_k^j,v_k^j\right\rangle_{L^2(\Omega_k\setminus\Omega_{k-1})}$$ $$=-{\rm Im}\,\left\langle \partial_\nu u_{m+1}^j,u_{m+1}^j\right\rangle_{L^2(\Gamma_{m+1})}=\lambda^j\,\left\langle a u_{m+1}^j,u_{m+1}^j\right\rangle_{L^2(\Gamma_{m+1})}.\eqno{(2.3)}$$
By (1.4) and (2.3) we conclude
$$a_0\lambda^j\left\|u_{m+1}^j|_{\Gamma_{m+1}}\right\|_{L^2(\Gamma_{m+1})}^2\le \gamma\lambda^j\|u^j\|_H^2+\gamma^{-1}\lambda^{-j}\|v^j\|_H^2,\eqno{(2.4)}$$
for every $\gamma>0$. On the other hand, applying (2.2) with $w=u^j$ yields
$$\|u^j\|_H^2\le C\lambda^{-2}\|v^j\|_H^2+C\left\|u_{m+1}^j|_{\Gamma_{m+1}}\right\|_{L^2(\Gamma_{m+1})}^2.\eqno{(2.5)}$$
Combining (2.4) and (2.5) and taking $\gamma$ small enough, independent of $\lambda$, we get
$$\|u^j\|_H\le C\lambda^{-j}\|v^j\|_H,\eqno{(2.6)}$$
which is equivalent to (1.7) for real $\lambda\gg 1$. Clearly, the case $-\lambda\gg 1$ can be treated in the same way.\\

{\it Proof of Theorem 2.1.} Given any $1\le k\le m$, define the resolvent ${\cal R}_k(\lambda)$ as $u={\cal R}_k(\lambda)v$, where $u=(u_1,...,u_{k+1})$, $v=(v_1,...,v_{k+1})$ satisfy the equation
$$
\left\{
\begin{array}{l}
 (\lambda^2+c_\ell^2\Delta)u_\ell=v_\ell\quad\mbox{in}\quad\Omega_\ell\setminus\Omega_{\ell-1},\,\ell=1,...,k,\\
 (\lambda^2+c_{k+1}^2\Delta)u_{k+1}=v_{k+1}\quad\mbox{in}\quad{\bf R}^n\setminus\Omega_k,\\
Bu_1=0\quad\mbox{on}\quad\Gamma_0,\\
u_\ell=u_{\ell+1},\,\partial_\nu u_\ell=\partial_\nu u_{\ell+1}
\quad\mbox{on}\quad \Gamma_\ell,\,\ell=1,...,k,\\
 u_{k+1} - \lambda-\mbox{outgoing}.
\end{array}
\right.
\eqno{(2.7)}
$$
Let us first see that Theorem 2.1 follows from the following

\begin{Theorem} Assume (1.3) and (1.6) fulfilled. Then, for every $1\le k\le m$ the cutoff resolvent $\chi{\cal R}_k(\lambda)\chi$ satisfies the estimate
$$\left\|\chi {\cal R}_k(\lambda)\chi\right\|_{L^2({\bf R}^n\setminus\Omega_k)\to L^2({\bf R}^n\setminus\Omega_k)}\le C|\lambda|^{-1}\quad\mbox{for}\quad \lambda\in{\bf R},\,|\lambda|\ge 1,\eqno{(2.8)}$$
where $\chi\in C_0^\infty({\bf R}^n)$, $\chi=1$ on $\Omega_k$.
\end{Theorem}

Choose a real-valued function $\rho\in C_0^\infty({\bf R})$, $0\le\rho\le 1$, $\rho(t)=1$ for $|t|\le\delta/2$, $\rho(t)=0$ for $|t|\ge\delta$,
$d\rho(t)/dt\le 0$ for $t\ge 0$, where $0<\delta\ll 1$ is a parameter. Given $x\in \Omega_{m+1}\setminus\Omega_m$, denote by $d(x)$ the distance between $x$ and $\Gamma_{m+1}$. Hence $\psi(x)=\rho(d(x))\in C^\infty(\Omega_{m+1})$, $\psi=1$ near $\Gamma_{m+1}$, $\psi=0$ on $\Omega_m$. The following estimate is proved in \cite{kn:CPV} (see Proposition 2.2) using in an essential way that the boundary $\Gamma_{m+1}$ is strictly concave viewed from the interior.

\begin{prop} There exist constants $C,\lambda_0,\delta_0>0$ so that if $0<\delta\le\delta_0$, $\lambda\ge\lambda_0$, we have the estimate
$$\|\psi u\|_{H^1(\Omega_{m+1}\setminus\Omega_m)}\le C\lambda^{-1}\|(\lambda^2+c_{m+1}^2\Delta) u\|_{L^2(\Omega_{m+1}\setminus\Omega_m)}$$ 
$$+C\left\|u|_{\Gamma_{m+1}}\right\|_{L^2(\Gamma_{m+1})}+C\lambda^{-1}\left\|\partial_\nu u|_{\Gamma_{m+1}}\right\|_{L^2(\Gamma_{m+1})}$$ $$
+ O_\delta(\lambda^{-1/2})\| u\|_{H^1(\Omega_{m+1}\setminus\Omega_m)},\quad\forall u\in H^2(\Omega_{m+1}\setminus\Omega_m),\eqno{(2.9)}$$
where the Sobolev space $H^1$ is equipped with the semi-classical norm with a small parameter $\lambda^{-1}$.
\end{prop}

Let $\chi\in C_0^\infty({\bf R}^n)$, $\chi=1$ on $\Omega_m$, supp$\,\chi\subset\Omega_{m+1}$. Clearly, the solution to (2.1) satisfies the equation
$$
\left\{
\begin{array}{l}
 (\lambda^2+c_k^2\Delta)w_k=v_k\quad\mbox{in}\quad\Omega_k\setminus\Omega_{k-1},\,k=1,...,m,\\
 (\lambda^2+c_{m+1}^2\Delta)\chi w_{m+1}=\chi v_{m+1}+c_{m+1}^2[\Delta,\chi] w_{m+1}\quad\mbox{in}\quad{\bf R}^n\setminus\Omega_m,\\
Bw_1=0\quad\mbox{on}\quad\Gamma_0,\\
w_k=w_{k+1},\,\partial_\nu w_k=\partial_\nu w_{k+1}
\quad\mbox{on}\quad \Gamma_k,\,k=1,...,m.
\end{array}
\right.
\eqno{(2.10)}
$$
Therefore, applying (2.8) with $k=m$ leads to the estimate
$$\sum_{k=1}^m\|w_k\|_{L^2(\Omega_k\setminus\Omega_{k-1})}+\|\chi w_{m+1}\|_{L^2({\bf R}^n\setminus\Omega_m)}$$ $$\le C\lambda^{-1}\sum_{k=1}^{m+1}\|v_k\|_{L^2(\Omega_k\setminus\Omega_{k-1})}+C\lambda^{-1}\|[\Delta,\chi] w_{m+1}\|_{L^2(\Omega_{m+1}\setminus\Omega_m)}.\eqno{(2.11)}$$
Choose $\chi$ so that $\psi=1$ on both supp$\,[\Delta,\chi]$ and supp$\,(1-\chi)|_{\Omega_{m+1}}$. Then (2.11) can be rewritten as follows
$$\sum_{k=1}^{m+1}\|w_k\|_{L^2(\Omega_k\setminus\Omega_{k-1})}\le C\lambda^{-1}\sum_{k=1}^{m+1}\|v_k\|_{L^2(\Omega_k\setminus\Omega_{k-1})}+C\|\psi w_{m+1}\|_{H^1(\Omega_{m+1}\setminus\Omega_m)},\eqno{(2.12)}$$
where again $H^1$ is equipped with the semiclassical norm. Using (2.9) with $u=w_{m+1}$ and combining with (2.12) lead to the estimate
$$\sum_{k=1}^{m+1}\|w_k\|_{L^2(\Omega_k\setminus\Omega_{k-1})}\le C\lambda^{-1}\sum_{k=1}^{m+1}\|v_k\|_{L^2(\Omega_k\setminus\Omega_{k-1})}+C\lambda^{-1/2}\| w_{m+1}\|_{H^1(\Omega_{m+1}\setminus\Omega_m)}$$ $$+C\left\|w_{m+1}|_{\Gamma_{m+1}}\right\|_{L^2(\Gamma_{m+1})}+C\lambda^{-1}\left\|\partial_\nu w_{m+1}|_{\Gamma_{m+1}}\right\|_{L^2(\Gamma_{m+1})}.\eqno{(2.13)}$$
On the other hand, by Green's formula we have
$$\sum_{k=1}^{m+1}{\rm Re}\,\left\langle w_k,v_k\right\rangle_{L^2(\Omega_k\setminus\Omega_{k-1},c_k^{-2}dx)}=\lambda^2\sum_{k=1}^{m+1}\|w_k\|^2_{L^2(\Omega_k\setminus\Omega_{k-1},c_k^{-2}dx)}$$ $$-\sum_{k=1}^{m+1}\|\nabla w_k\|^2_{L^2(\Omega_k\setminus\Omega_{k-1})}-{\rm Re}\,\left\langle \partial_\nu w_{m+1},w_{m+1}\right\rangle_{L^2(\Gamma_{m+1})},$$
which in turn implies
$$\sum_{k=1}^{m+1}\|\nabla w_k\|_{L^2(\Omega_k\setminus\Omega_{k-1})}\le C\lambda^{-1}\sum_{k=1}^{m+1}\|v_k\|_{L^2(\Omega_k\setminus\Omega_{k-1})}+C\sum_{k=1}^{m+1}\| w_k\|_{L^2(\Omega_k\setminus\Omega_{k-1})}$$
  $$+C\lambda^{-1/2}\left(\left\|w_{m+1}|_{\Gamma_{m+1}}\right\|_{L^2(\Gamma_{m+1})}+\left\|\lambda^{-1}\partial_\nu w_{m+1}|_{\Gamma_{m+1}}\right\|_{L^2(\Gamma_{m+1})}\right).\eqno{(2.14)}$$
Combining (2.13) and (2.14) and taking $\lambda$ large enough, we conclude that the second term in the right-hand side of (2.13) can be absorbed, thus obtaining (2.2).
\eproof

{\it Proof of Theorem 2.2.} Since (2.8) holds true for $k=0$ in view of the assumption (1.6), one needs to show that (2.8) with $k-1$ implies (2.8) with $k$. This, however, is proved in \cite{kn:CPV} (see Theorem 1.1; see also Section 4 below).
\eproof

The fact that (1.7) implies (1.8) and (1.9) is more or less well known. In what follows we will sketch the main points. Define the operator $A_0$ on the Hilbert space ${\cal H}_0=H$ as follows
$$A_0u=(-c_1^2\Delta u_1,...,-c_{m+1}^2\Delta u_{m+1}),\quad u=(u_1,...,u_{m+1}),$$
with domain of definition
$${\cal D}(A_0)=\left\{u\in H: A_0u\in H, Bu_1|_{\Gamma_0}=0, u_k|_{\Gamma_k}=u_{k+1}|_{\Gamma_k}, \partial_\nu u_k|_{\Gamma_k}=\partial_\nu u_{k+1}|_{\Gamma_k}, k=1,...,m,\right.$$ $$\left. \partial_\nu u_{m+1}|_{\Gamma_{m+1}}=-ia(x)u_{m+1}|_{\Gamma_{m+1}}\right\}.$$
By Green's formula we have
$${\rm Im}\,\left\langle A_0u,u\right\rangle_H=-{\rm Im}\,\left\langle\partial_\nu u_{m+1},u_{m+1}\right\rangle_{L^2(\Gamma_{m+1})}=
\left\langle a u_{m+1},u_{m+1}\right\rangle_{L^2(\Gamma_{m+1})}\ge 0,$$
which in turn implies that $A_0$ is a generator of a semi-group $e^{itA_0}$. Then the solutions to (1.1) can be expressed by the formula
$$u^0(t)=e^{itA_0}u^0(0),\quad t\ge 0.$$
It follows from \cite{kn:B} that, under the assumption (1.4), $A_0$ has no eigenvalues on the real axis. Moreover, applying (1.7) with $j=0$ and $z=\lambda^2$ yields that the resolvent $(A_0-z)^{-1}$ is analytic in a strip $|{\rm Im}\,z|\le\gamma_0$, $\gamma_0>0$, and satisfies in this region the bound
$$\left\|(A_0-z)^{-1}\right\|_{{\cal H}_0\to {\cal H}_0}\le Const,$$
which in turn implies
$$\left\|e^{itA_0}\right\|_{{\cal H}_0\to {\cal H}_0}\le \widetilde C e^{-Ct},\quad t>0,\eqno{(2.15)}$$
with constants $\widetilde C,C>0$ independent of $t$. Clearly, (2.15) is equivalent to (1.8).

We would like to treat the equation (1.2) in a similar way. To this end, introduce the Hilbert space ${\cal H}_1=\dot H^1_B\oplus H$, where
$$\dot H^1_B=\dot H^1_B(\Omega_1\setminus\Omega_0)\oplus\oplus_{k=2}^{m+1}\dot H^1(\Omega_k\setminus\Omega_{k-1}),$$
$$\dot H^1(\Omega_k\setminus\Omega_{k-1})=\left\{u:\,\int_{\Omega_k\setminus\Omega_{k-1}}|\nabla u|^2dx<+\infty\right\},\quad 2\le k\le m+1,$$
$$\dot H^1_B(\Omega_1\setminus\Omega_0)=\left\{u:\,\int_{\Omega_1\setminus\Omega_0}|\nabla u|^2dx<+\infty\right\},\quad \mbox{if}\quad B=\partial_\nu,$$
$$\dot H^1_B(\Omega_1\setminus\Omega_0)=\left\{u:\,\int_{\Omega_1\setminus\Omega_0}|\nabla u|^2dx<+\infty,\,u|_{\Gamma_0}=0\right\},\quad \mbox{if}\quad B=Id.$$
On ${\cal H}_1$ define the operator $A_1$ as follows
$$A_1= -i\left( \begin{array}{ll} \hskip 0.5cm 0& Id\\ c^2(x) \Delta &0
\end{array}
\right),$$ where $$c^2(x) \Delta u:=(c_1^2\Delta u_1,...,c_{m+1}^2\Delta u_{m+1}),\quad u=(u_1,...,u_{m+1}),$$
with domain of definition 
$${\cal D}(A_1)=\left\{(u,v)\in {\cal H}_1: v\in \dot H^1_B, c^2(x) \Delta u\in H,  Bu_1|_{\Gamma_0}=0, u_k|_{\Gamma_k}=u_{k+1}|_{\Gamma_k},  \right.$$ $$\left. \partial_\nu u_k|_{\Gamma_k}=\partial_\nu u_{k+1}|_{\Gamma_k},k=1,...,m,\,\partial_\nu u_{m+1}|_{\Gamma_{m+1}}=-a(x)v_{m+1}|_{\Gamma_{m+1}}\right\}.$$
By Green's formula we have
$${\rm Im}\,\left\langle A_1\left(
\begin{array}{ll} u\\ v
\end{array}\right),\left(
\begin{array}{ll} u\\ v
\end{array}\right)\right\rangle_{{\cal H}_1}=-{\rm Re}\,\left\langle\left(
\begin{array}{ll} v\\ c^2(x)\Delta u
\end{array}\right),\left(
\begin{array}{ll} u\\ v
\end{array}\right)\right\rangle_{{\cal H}_1}$$ $$=-{\rm Re}\,\left\langle\partial_\nu u_{m+1},v_{m+1}\right\rangle_{L^2(\Gamma_{m+1})}=
\left\langle a v_{m+1},v_{m+1}\right\rangle_{L^2(\Gamma_{m+1})}\ge 0,$$
which in turn implies that $A_1$ is a generator of a semi-group $e^{itA_1}$. Then the solutions to (1.2) can be expressed by the formula
$$\left(
\begin{array}{ll} u^1(t)\\ \partial_tu^1(t)
\end{array}\right)=e^{itA_1}\left(\begin{array}{ll} u^1(0)\\ \partial_tu^1(0)
\end{array}\right),\quad t\ge 0.$$
It follows from \cite{kn:B} that, under the assumption (1.4), $A_1$ has no eigenvalues on the real axis. Moreover, applying (1.7) with $j=1$ and $z=\lambda$ yields that the resolvent $(A_1-z)^{-1}$ is analytic in a strip $|{\rm Im}\,z|\le\gamma_1$, $\gamma_1>0$, and satisfies in this region the bound
$$\left\|(A_1-z)^{-1}\right\|_{{\cal H}_1\to {\cal H}_1}\le Const,$$
which in turn implies
$$\left\|e^{itA_1}\right\|_{{\cal H}_1\to {\cal H}_1}\le \widetilde C e^{-Ct},\quad t>0,\eqno{(2.16)}$$
with constants $\widetilde C,C>0$ independent of $t$. It is easy to see that (2.16) is equivalent to (1.9).

Introduce the Hilbert space ${\cal H}=\dot H^1_{B,sc}\oplus H_{sc}$, where
$$H_{sc}:=\oplus_{k=1}^{m}L^2\left(\Omega_{k}\setminus\Omega_{k-1},c_k^{-2}dx\right)\oplus L^2\left({\bf R}^n\setminus\Omega_m,c_{m+1}^{-2}dx\right),$$
$$\dot H^1_{B,sc}=\dot H^1_B(\Omega_1\setminus\Omega_0)\oplus\oplus_{k=2}^{m}\dot H^1(\Omega_k\setminus\Omega_{k-1})\oplus\dot H^1\left({\bf R}^n\setminus\Omega_m\right),$$
$$\dot H^1({\bf R}^n\setminus\Omega_m)=\left\{u:\,\int_{{\bf R}^n\setminus\Omega_m}|\nabla u|^2dx<+\infty\right\}.$$
On ${\cal H}$ define the operator $A$ as follows
$$A= -i\left( \begin{array}{ll} \hskip 0.5cm 0& Id\\ c^2(x) \Delta &0
\end{array}
\right), $$
with domain of definition 
$${\cal D}(A)=\left\{(u,v)\in {\cal H}: v\in \dot H^1_{B,sc}, c^2(x) \Delta u\in H_{sc},  Bu_1|_{\Gamma_0}=0, u_k|_{\Gamma_k}=u_{k+1}|_{\Gamma_k},  \right.$$ $$\left. \partial_\nu u_k|_{\Gamma_k}=\partial_\nu u_{k+1}|_{\Gamma_k},k=1,...,m\right\}.$$
By Green's formula we have
$${\rm Im}\,\left\langle A\left(
\begin{array}{ll} u\\ v
\end{array}\right),\left(
\begin{array}{ll} u\\ v
\end{array}\right)\right\rangle_{{\cal H}}=0,$$
which in turn implies that $A$ is a generator of a group $e^{itA}$. Then the solutions to (1.17) can be expressed by the formula
$$\left(
\begin{array}{ll} u(t)\\ \partial_tu(t)
\end{array}\right)=e^{itA}\left(\begin{array}{ll} u(0)\\ \partial_tu(0)
\end{array}\right),\quad t\ge 0.$$
As in \cite{kn:V}, it follows from (2.8) applied with $k=m$ and $z=\lambda$ that the cutoff resolvent $\chi(A-z)^{-1}\chi$ is analytic in a strip $|{\rm Im}\,z|\le\gamma$, $\gamma>0$, and satisfies in this region the bound
$$\left\|\chi(A-z)^{-1}\chi\right\|_{{\cal H}\to {\cal H}}\le Const,$$
where $\chi\in C_0^\infty({\bf R}^n)$, $\chi=1$ on $\Omega_m$. 
This in turn implies (see \cite{kn:V}, \cite{kn:K})
$$\left\|\chi e^{itA}\chi\right\|_{{\cal H}\to {\cal H}}\le C_\chi p_0(t),\quad t>0,\eqno{(2.17)}$$
with a constant $C_\chi>0$ independent of $t$. It is easy to see that (2.17) is equivalent to (1.18).

\section{The case $\Omega_0$ trapping}

As in the previous section, Theorem 1.3 follows from the following

\begin{Theorem} Assume (1.3) and (1.12) fulfilled. Then, there exist constants $C,\lambda_0>0$ so that for $\lambda\ge\lambda_0$ the solution to (2.1) satisfies the estimate
$$(\log\lambda)^{-2^m}\|w\|_H\le C\lambda^{-1}\|v\|_H+C\left\|w_{m+1}|_{\Gamma_{m+1}}\right\|_{L^2(\Gamma_{m+1})}+C\lambda^{-1}\left\|\partial_\nu w_{m+1}|_{\Gamma_{m+1}}\right\|_{L^2(\Gamma_{m+1})}.\eqno{(3.1)}$$
\end{Theorem}

Moreover, proceeding as in Section 2 it is easy to see that Theorem 3.1 follows from the following theorem the proof of which will be given in the next section.

\begin{Theorem} Assume (1.3) and (1.12) fulfilled. Then, for every $0\le k\le m$ the cutoff resolvent $\chi{\cal R}_k(\lambda)\chi$ is analytic in $\{\lambda\in{\bf C}:|{\rm Im}\,\lambda|\le C_1(\log|\lambda|)^{-2^k},\,|\lambda|\ge C_2\}$ and satisfies in this region the estimate
$$\left\|\chi {\cal R}_k(\lambda)\chi\right\|_{L^2({\bf R}^n\setminus\Omega_k)\to L^2({\bf R}^n\setminus\Omega_k)}\le C|\lambda|^{-1}(\log|\lambda|)^{2^k},\eqno{(3.2)}$$
where $C,C_1$ and $C_2$ are positive constants.
\end{Theorem}

\noindent
{\bf Remark.} It is natural to expect that (1.12) implies that all cutoff resolvents $\chi{\cal R}_k(\lambda)\chi$, $k=1,...,m$, are analytic in
some strip $\{|{\rm Im}\,\lambda|\le C_1,|\lambda|\ge C_2\}$, $C_1,C_2>0$. However, this remains a difficult open problem. Note that large free of resonances regions far from the real axis are obtained in \cite{kn:CPV2} under some natural assumptions.

To prove Corollary 1.4 observe first that (1.13) is equivalent to the estimate (with $j=0,1$)
$$\left\|(A_j-z)^{-1}\right\|_{{\cal H}_j\to {\cal H}_j}\le C\left(\log|z|\right)^{2^{m+1}}\quad{\rm for}\,\, z\in{\bf R},\,|z|\ge C',\eqno{(3.3)}$$
with some constants $C>0$, $C'>2$ independent of $z$. Clearly, (3.3) implies that $(A_j-z)^{-1}$ is analytic in
$$\Lambda=\left\{z\in {\bf C}: |{\rm Im}\,z|\le C_1\left(\log|z|\right)^{-2^{m+1}},\,|z|\ge C_2\right\}$$
and satisfies in this region the bound (3.3). Therefore, using the fact that the operators $A_j$ are elliptic together with a standard interpolation argument, we conclude that
$$\left\|(A_j-z)^{-1}\right\|_{{\cal H}_j^\varepsilon\to {\cal H}_j}\le C_\varepsilon\quad{\rm for}\,\, z\in\Lambda,\eqno{(3.4)}$$
for every $\varepsilon>0$ with a constant $C_\varepsilon>0$ independent of $z$, where ${\cal H}_0^\varepsilon:= H^\varepsilon$, while the norm
$\|\cdot\|_{{\cal H}_1^\varepsilon}$ is defined by replacing in the definition of ${\cal H}_1$ all norms $L^2$ by the Sobolev norms $H^\varepsilon$.
On the other hand, proceeding as in \cite{kn:L} one can show that (3.4) implies
$$\left\|e^{itA_j}\right\|_{{\cal H}_j^\varepsilon\to {\cal H}_j}\le \widetilde C_\varepsilon\exp\left(-C\varepsilon t^{\alpha_{m+1}}\right),\quad t>0,\eqno{(3.5)}$$
for $0<\varepsilon\le\varepsilon_0$, with constants $C,\varepsilon_0>0$ independent of $t$ and $\varepsilon$, $\widetilde C_\varepsilon>0$ independent of $t$. Clearly, (3.5) is equivalent to (1.14) and (1.15), respectively.

Similarly, the estimate (3.2) with $k=m$ implies that the cutoff resolvent $\chi(A-z)^{-1}\chi$ is analytic in $\{z\in{\bf C}: |{\rm Im}\,z|\le C_1(\log|z|)^{-2^m},\,|z|\ge C_{2}\}$ and satisfies in this region the estimate
$$\left\|\chi(A-z)^{-1}\chi\right\|_{{\cal H}^\varepsilon\to {\cal H}}\le C_\varepsilon,\eqno{(3.6)}$$
where ${\cal H}^\varepsilon$ is defined as ${\cal H}_1^\varepsilon$ above. On the other hand, as in \cite{kn:PV1} one can show that (3.6) implies
$$\left\|\chi e^{itA}\chi\right\|_{{\cal H}^\varepsilon\to {\cal H}}\le C_{\chi,\varepsilon} p_\varepsilon(t),\quad t>0,\eqno{(3.7)}$$
with a constant $C_{\chi,\varepsilon}>0$ independent of $t$. It is easy to see that (3.7) is equivalent to (1.19).

\section{Proof of Theorem 3.2}

We will prove (3.2) by induction in $k$. Let us first see that the assumption (1.12) implies (3.2) with $k=0$. This is essentially proved in \cite{kn:Bu} (see Proposition 4.4 and Lemma 4.7). The idea is to apply the Phragm\`en-Lindel\"of principle to the operator-valued function
$$g(\lambda)=\frac{\lambda e^{iN\lambda\log\lambda}}{\log\lambda}\chi{\cal R}_0(\lambda)\chi,\quad {\rm Re}\,\lambda\ge C_2,$$
where $\log\lambda=\log|\lambda|+i\arg\lambda$ and $N>0$ is a constant big enough. It is well known that the outgoing resolvent satisfies the bound
$$\left\|{\cal R}_0(\lambda)\right\|_{L^2({\bf R}^n\setminus\Omega_0)\to L^2({\bf R}^n\setminus\Omega_0)}\le\frac{1}{|\lambda||{\rm Im}\,\lambda|}\quad\mbox{for}\quad {\rm Im}\,\lambda<0.\eqno{(4.1)}$$
Hence, on ${\rm Im}\,\lambda=-(N\log|\lambda|)^{-1}$, ${\rm Re}\,\lambda\ge C_2$, we have the bound
$$\left\|g(\lambda)\right\|_{L^2\to L^2}\le \frac{Ce^{-N{\rm Im}\,(\lambda\log\lambda)}}{|{\rm Im}\,\lambda|\log|\lambda|}\le \frac{Ce^{N|{\rm Im}\,\lambda|\log|\lambda|}}{|{\rm Im}\,\lambda|\log|\lambda|}\le Const.\eqno{(4.2)}$$
On the other hand, by (1.12), on ${\rm Im}\,\lambda=C_1>0$, ${\rm Re}\,\lambda\ge C_2$, we have the bound
$$\left\|g(\lambda)\right\|_{L^2\to L^2}\le C|\lambda|^{p+1}e^{-N{\rm Im}\,(\lambda\log\lambda)}\le Ce^{(p+1-N{\rm Im}\,\lambda)\log|\lambda|}\le Const,\eqno{(4.3)}$$
if we choose $N=(p+1)/C_1$. By the Phragm\`en-Lindel\"of principle, we conclude from (4.2) and (4.3) that the function $g(\lambda)$ satisfies the bound
$$\left\|g(\lambda)\right\|_{L^2\to L^2}\le Const,\eqno{(4.4)}$$
in $-(N\log|\lambda|)^{-1}\le {\rm Im}\,\lambda\le C_1$, ${\rm Re}\,\lambda\ge C_2$. It follows from (4.4) that for $-(N\log|\lambda|)^{-1}\le {\rm Im}\,\lambda\le \varepsilon/2N$, ${\rm Re}\,\lambda\ge C_2$, $0<\varepsilon\ll 1$, we have
$$\left\|\lambda\chi{\cal R}_0(\lambda)\chi\right\|_{L^2\to L^2}\le C\log|\lambda|e^{N{\rm Im}\,(\lambda\log\lambda)}\le C\log|\lambda|e^{\frac{\varepsilon \log|\lambda|}{2}}\le C\varepsilon^{-1}|\lambda|^\varepsilon,\eqno{(4.5)}$$
with a constant $C>0$ independent of $\lambda$ and $\varepsilon$. On the other hand, for $-\varepsilon/2N\le {\rm Im}\,\lambda\le -(N\log|\lambda|)^{-1}$ the estimate (4.5) follows from (4.1). Thus we conclude that (4.5) holds for $|{\rm Im}\,\lambda|\le \varepsilon/2N$, ${\rm Re}\,\lambda\ge C_2$. Clearly, the case ${\rm Re}\,\lambda\le -C_2$ can be treated similarly. Taking $\varepsilon$ such that $|\lambda|^\varepsilon=2$, we obtain (3.2) with $k=0$.

Thus, to prove Theorem 3.2 it suffices to show that (3.2) with $k-1$, $1\le k\le m$, implies (3.2) with $k$. 
Let $w=(w_1,...,w_k)$, $v=(v_1,...,v_k)$ satisfy the equation
$$
\left\{
\begin{array}{l}
 (\lambda^2+c_\ell^2\Delta)w_\ell=v_\ell\quad\mbox{in}\quad\Omega_\ell\setminus\Omega_{\ell-1},\,\ell=1,...,k,\\
Bw_1=0\quad\mbox{on}\quad\Gamma_0,\\
w_\ell=w_{\ell+1},\,\partial_\nu w_\ell=\partial_\nu w_{\ell+1}
\quad\mbox{on}\quad \Gamma_\ell,\,\ell=1,...,k-1.
\end{array}
\right.
\eqno{(4.6)}
$$
We need the following extension of Theorem 2.1.

\begin{Theorem} Assumed (3.2) fulfilled with $k-1$. Then, there exist constants $C,\lambda_0>0$ so that for  $\lambda\ge\lambda_0$ 
 the solution to (4.6) satisfies the estimate
$$(\log\lambda)^{-2^{k-1}}\|w\|_{H_k}\le C\lambda^{-1}\|v\|_{H_k}+C\left\|w_k|_{\Gamma_k}\right\|_{L^2(\Gamma_k)}+C\lambda^{-1}\left\|\partial_\nu w_k|_{\Gamma_k}\right\|_{L^2(\Gamma_k)},\eqno{(4.7)}$$
where $H_k:=\oplus_{\ell=1}^kL^2\left(\Omega_\ell\setminus\Omega_{\ell-1},c_\ell^{-2}dx\right)$.
\end{Theorem}

{\it Proof.} Let $\chi\in C_0^\infty({\bf R}^n)$, $\chi=1$ on $\Omega_{k-1}$, supp$\,\chi\subset\Omega_k$, such that $\psi=1$ on supp$\,[\Delta,\chi]$ and supp$\,(1-\chi)|_{\Omega_k}$. We have
$$\chi w=\chi_1{\cal R}_{k-1}(\lambda)\chi_1\left(\chi v+[\Delta,\chi]w_k\right),\eqno{(4.8)}$$
where $\chi_1=1$ on supp$\,\chi$, supp$\,\chi_1\subset\Omega_k$. By (3.2) with $k-1$ and (4.8) we conclude
$$(\log\lambda)^{-2^{k-1}}\|w\|_{H_k}\le(\log\lambda)^{-2^{k-1}}\left(\|\chi w\|_{H_k}+\|\psi w_k\|_{L^2(\Omega_k\setminus\Omega_{k-1})}\right)$$ $$\le C\lambda^{-1}\|v\|_{H_k}+C\|\psi w_k\|_{H^1(\Omega_k\setminus\Omega_{k-1})},\eqno{(4.9)}$$
where $H^1$ is equipped with the semiclassical norm. By (2.9) and (4.9),
$$(\log\lambda)^{-2^{k-1}}\|w\|_{H_k}\le C\lambda^{-1}\|v\|_{H_k}+C\lambda^{-1/2}\| w_k\|_{H^1(\Omega_k\setminus\Omega_{k-1})}$$ $$+C\left\|w_k|_{\Gamma_k}\right\|_{L^2(\Gamma_k)}+C\lambda^{-1}\left\|\partial_\nu w_k|_{\Gamma_k}\right\|_{L^2(\Gamma_k)}.\eqno{(4.10)}$$
On the other hand, we have an analogue of (2.14) with $m+1$ replaced by $k$, which together with (4.10) yield
$$(\log\lambda)^{-2^{k-1}}\left(\|w\|_{H_k}+\| w_k\|_{H^1(\Omega_k\setminus\Omega_{k-1})}\right)\le C\lambda^{-1}\|v\|_{H_k}+C\lambda^{-1/2}\| w_k\|_{H^1(\Omega_k\setminus\Omega_{k-1})}$$ $$+C\left\|w_k|_{\Gamma_k}\right\|_{L^2(\Gamma_k)}+C\lambda^{-1}\left\|\partial_\nu w_k|_{\Gamma_k}\right\|_{L^2(\Gamma_k)}.\eqno{(4.11)}$$
Clearly, we can absorb the second term in the right-hand side of (4.11) by taking $\lambda$ big enough, thus obtaining (4.7).
\eproof

Note that it suffices to prove (3.2) for $\lambda\in{\bf R}$, $|\lambda|\gg 1$, only (see \cite{kn:V}). Without loss of generality we may suppose $\lambda>0$. 
Let $u=(u_1,...,u_{k+1})$, $v=(v_1,...,v_{k+1})$ satisfy the equation (2.7) with supp$\,v_{k+1}\subset K$, where $K\subset{\bf R}^n\setminus\Omega_k$ is a compact. Set $f_k=u_{k+1}|_{\Gamma_k}=u_{k}|_{\Gamma_k}$. Define the outgoing Neumann operator, $N_k(\lambda)$, for the exterior problem in ${\bf R}^n\setminus\Omega_k$ as follows
$$N_k(\lambda)f=\lambda^{-1}\partial_{\nu'}U_k(\lambda)f|_{\Gamma_k},$$
where $\nu'$ is the outer unit normal to $\Gamma_k$, and $U_k(\lambda)$ solves the equation
$$
\left\{
\begin{array}{l}
 (\lambda^2+c_{k+1}^2\Delta)U_k(\lambda)f=0\quad\mbox{in}\quad{\bf R}^n\setminus\Omega_k,\\
U_k(\lambda)f=f\quad\mbox{on}\quad\Gamma_k,\\
U_k(\lambda)f - \lambda - \mbox{outgoing}.
\end{array}
\right.
\eqno{(4.12)}
$$
Define also the operator $G_k(\lambda)$ via the equation
$$
\left\{
\begin{array}{l}
 (\lambda^2+c_{k+1}^2\Delta)G_k(\lambda)f=f\quad\mbox{in}\quad{\bf R}^n\setminus\Omega_k,\\
G_k(\lambda)f=0\quad\mbox{on}\quad\Gamma_k,\\
G_k(\lambda)f - \lambda - \mbox{outgoing}.
\end{array}
\right.
\eqno{(4.13)}
$$
Set $\widetilde u_{k+1}=u_{k+1}-G_k(\lambda)v_{k+1}$. Then, the equation (2.7) can be rewritten as follows
$$
\left\{
\begin{array}{l}
 (\lambda^2+c_\ell^2\Delta)u_\ell=v_\ell\quad\mbox{in}\quad\Omega_\ell\setminus\Omega_{\ell-1},\,\ell=1,...,k,\\
 (\lambda^2+c_{k+1}^2\Delta)\widetilde u_{k+1}=0\quad\mbox{in}\quad{\bf R}^n\setminus\Omega_k,\\
Bu_1=0\quad\mbox{on}\quad\Gamma_0,\\
u_\ell=u_{\ell+1},\,\partial_\nu u_\ell=\partial_\nu u_{\ell+1}
\quad\mbox{on}\quad \Gamma_\ell,\,\ell=1,...,k-1,\\
\widetilde u_{k+1}=u_k,\,\partial_{\nu'}\widetilde u_{k+1}=-\partial_\nu u_k+\lambda h_k,\quad\mbox{on}\quad\Gamma_k,\\
\widetilde u_{k+1} - \lambda-\mbox{outgoing},
\end{array}
\right.
\eqno{(4.14)}
$$
where $h_k=\lambda^{-1}\partial_{\nu'}G_k(\lambda)v_{k+1}|_{\Gamma_k}$, and we have used that $\nu'=-\nu$. Hence $\widetilde u_{k+1}=U_k(\lambda)f_k$, and (4.14) implies
$$
\left\{
\begin{array}{l}
 (\lambda^2+c_\ell^2\Delta)u_\ell=v_\ell\quad\mbox{in}\quad\Omega_\ell\setminus\Omega_{\ell-1},\,\ell=1,...,k,\\
Bu_1=0\quad\mbox{on}\quad\Gamma_0,\\
u_\ell=u_{\ell+1},\,\partial_\nu u_\ell=\partial_\nu u_{\ell+1}
\quad\mbox{on}\quad \Gamma_\ell,\,\ell=1,...,k-1,\\
 u_{k}=f_k,\,\lambda^{-1}\partial_{\nu} u_{k}=-N_k(\lambda)f_k+h_k,\quad\mbox{on}\quad\Gamma_k.\\
\end{array}
\right.
\eqno{(4.15)}
$$
The fact that $\Omega_k$ is strictly convex implies the bounds (see Theorem 3.1 of \cite{kn:CPV}):
$$\|h_k\|_{L^2(\Gamma_k)}\le C_K\lambda^{-1}\|v_{k+1}\|_{L^2({\bf R}^n\setminus\Omega_k)},\eqno{(4.16)}$$
$$\|u_{k+1}\|_{L^2(K)}\le \|U_{k}(\lambda)f_k\|_{L^2(K)}+\|G_{k}(\lambda)v_{k+1}\|_{L^2(K)}$$ $$\le C_K\|f_k\|_{H^1(\Gamma_k)}+C_K\lambda^{-1}\|v_{k+1}\|_{L^2({\bf R}^n\setminus\Omega_k)}.\eqno{(4.17)}$$
Hereafter all Sobolev spaces $H^1$ will be equipped with the semi-classical norm. Applying Green's formula to the solutions of (4.15) leads to the identity
$$-\lambda{\rm Im}\,\left\langle N_k(\lambda)f_k,f_k\right\rangle_{L^2(\Gamma_k)}+\lambda{\rm Im}\,\left\langle h_k,f_k\right\rangle_{L^2(\Gamma_k)}={\rm Im}\,\left\langle \partial_\nu u_k|_{\Gamma_k},u_k|_{\Gamma_k}\right\rangle_{L^2(\Gamma_k)}$$ 
$$=-\sum_{\ell=1}^k{\rm Im}\,\left\langle u_\ell,v_\ell\right\rangle_{L^2(\Omega_\ell\setminus\Omega_{\ell-1},c_\ell^{-2}dx)}.\eqno{(4.18)}$$
Hence, $\forall\beta>0$, we have
$$-{\rm Im}\,\left\langle N_k(\lambda)f_k,f_k\right\rangle_{L^2(\Gamma_k)}\le \beta^2\|f_k\|_{L^2(\Gamma_k)}^2+\beta^{-2}\|h_k\|_{L^2(\Gamma_k)}^2+\beta^2\|u\|_{H_k}^2+\beta^{-2}\lambda^{-2}\|v\|_{H_k}^2.\eqno{(4.19)}$$
Since $\Omega_k$ is strictly convex, the Neumann operator satisfies the bound (e.g. see Corollary 3.3 of \cite{kn:CPV})
$$\left\| N_k(\lambda)f_k\right\|_{L^2(\Gamma_k)}\le C\|f_k\|_{H^1(\Gamma_k)}.\eqno{(4.20)}$$
Applying Theorem 4.1 with $w=(u_1,...,u_k)$ and using (4.20), we get
$$(\log\lambda)^{-2^{k-1}}\|u\|_{H_k}\le C\lambda^{-1}\|v\|_{H_k}+C\|f_k\|_{H^1(\Gamma_k)}.\eqno{(4.21)}$$
Choose a function $\eta_k\in C^\infty(T^*\Gamma_k)$ such that $\eta_k=1$ on $\{\zeta\in T^*\Gamma_k:\|\zeta\|\le c_k^{-1}+\epsilon\}$,
$\eta_k=0$ on $\{\zeta\in T^*\Gamma_k:\|\zeta\|\le c_{k+1}^{-1}-\epsilon\}$, $0<\epsilon\ll 1$, which is possible in view of (1.3).
Recall that $\|\zeta\|^2$ is the principal symbol of the (positive) Laplace-Beltrami operator on $\Gamma_k$ evaluated at $\zeta$. 
We will denote by ${\rm Op}_\lambda(\eta_k)$ the $\lambda-\Psi$DO on $\Gamma_k$ with symbol $\eta_k$. Since $\Omega_k$ is strictly convex and $\eta_k$ is supported in the hyperbolic region for the corresponding exterior boundary value problem, it is well known that $N_k(\lambda){\rm Op}_\lambda(\eta_k)$ is a $\lambda-\Psi$DO with principal symbol $-i\eta_k(\zeta)\sqrt{c_{k+1}^{-2}-\|\zeta\|^2}$ (e.g. see the appendix of \cite{kn:G}). This together with (4.20) and G\"arding's inequality imply immediately the following

\begin{lemma} There exist constants $C_1,C_2>0$ such that we have
$$-{\rm Im}\,\left\langle N_k(\lambda)f_k,f_k\right\rangle_{L^2(\Gamma_k)}\ge C_1\|f_k\|_{L^2(\Gamma_k)}^2-C_2\|{\rm Op}_\lambda(1-\eta_k)f_k\|_{H^1(\Gamma_k)}^2.\eqno{(4.22)}$$
\end{lemma}

By (4.16), (4.17), (4.19), (4.21), (4.22), taking $\beta=\beta'(\log\lambda)^{-2^{k-1}}$ with $\beta'>0$ small enough independent of $\lambda$, we conclude
$$(\log\lambda)^{-2^{k-1}}\left(\|u\|_{H_k}+\|u_{k+1}\|_{L^2(K)}+\|f_k\|_{L^2(\Gamma_k)}\right)$$ $$\le C\lambda^{-1}(\log\lambda)^{2^{k-1}}\left(\|v\|_{H_k}+\|v_{k+1}\|_{L^2({\bf R}^n\setminus\Omega_k)}\right)+C\|{\rm Op}_\lambda(1-\eta_k)f_k\|_{H^1(\Gamma_k)}.\eqno{(4.23)}$$
On the other hand, the fact that $1-\eta_k$ is supported in the elliptic region for the corresponding interior boundary value problem implies the following

\begin{prop} There exist constants $C,\lambda_0>0$ so that for $\lambda\ge\lambda_0$ we have
$$\|{\rm Op}_\lambda(1-\eta_k)f_k\|_{H^1(\Gamma_k)}\le C\lambda^{-3/2}\|v_k\|_{L^2(\Omega_k\setminus\Omega_{k-1})}+C\lambda^{-1/2}\|u_k\|_{L^2(\Omega_k\setminus\Omega_{k-1})}$$ $$
+C\lambda^{-1/3}\|f_k\|_{L^2(\Gamma_k)}+C\|h_k\|_{L^2(\Gamma_k)}.\eqno{(4.24)}$$
\end{prop}

{\it Proof.} Choose a smooth function $\psi$ such that $\psi=1$ on $\{x: {\rm dist}(x,\Gamma_k)\le\delta\}$, $\psi=0$ outside $\{x: {\rm dist}(x,\Gamma_k)\le 2\delta\}$, where $\delta>0$ is a small parameter independent of $\lambda$. Set $\varphi(\zeta)=(1-\eta_k(\zeta))\langle\zeta\rangle$, $\zeta\in T^*\Gamma_k$, $w_k=\psi{\rm Op}_\lambda(\varphi)u_k$. Clearly,
$g_k:=w_k|_{\Gamma_k}={\rm Op}_\lambda(\varphi)f_k$, 
$$\lambda^{-1}\partial_\nu w_k|_{\Gamma_k}=\lambda^{-1}{\rm Op}_\lambda(\varphi)\partial_\nu u_k|_{\Gamma_k}=-{\rm Op}_\lambda(\varphi)N_k(\lambda)f_k+{\rm Op}_\lambda(\varphi)h_k$$
$$=-N_k(\lambda)g_k+[{\rm Op}_\lambda(\varphi),N_k(\lambda)]f_k+{\rm Op}_\lambda(\varphi)h_k.$$
By Green's formula we have
$$\lambda M+\lambda^{-1}{\rm Re}\,\left\langle\left(c_k^2\Delta+\lambda^{2}\right)w_k,w_k\right\rangle_{L^2(\Omega_k\setminus\Omega_{k-1},c_k^{-2}dx)}=
-{\rm Re}\,\left\langle\lambda^{-1}\partial_\nu w_k|_{\Gamma_k},w_k|_{\Gamma_k}\right\rangle_{L^2(\Gamma_k)}$$
$$={\rm Re}\,\left\langle N_k(\lambda)g_k,g_k\right\rangle_{L^2(\Gamma_k)}-{\rm Re}\,\left\langle [{\rm Op}_\lambda(\varphi),N_k(\lambda)]f_k,g_k\right\rangle_{L^2(\Gamma_k)}-{\rm Re}\,\left\langle {\rm Op}_\lambda(\varphi)h_k,g_k\right\rangle_{L^2(\Gamma_k)},\eqno{(4.25)}$$
where
$$M=\left\|\lambda^{-1}\nabla w_k\right\|_{L^2(\Omega_k\setminus\Omega_{k-1})}^2-c_k^{-2}\left\| w_k\right\|_{L^2(\Omega_k\setminus\Omega_{k-1})}^2.$$
Let us see that
$$\|g_k\|_{L^2(\Gamma_k)}^2\le C\lambda M,\quad C>0.\eqno{(4.26)}$$
Denote by $x_n>0$ the normal coordinate to $\Gamma_k$, i.e. given $x\in \Omega_k$, we have $x_n={\rm dist}(x,\Gamma_k)$. Given $0<x_n\le 2\delta\ll 1$, set $\Gamma_k(x_n)=\{x\in \Omega_k:{\rm dist}(x,\Gamma_k)=x_n\}$. Clearly, $M$ can be written in the form
$$M=\left\|\lambda^{-1}\partial_{x_n}w_k\right\|_{L^2}^2+\left\langle\left(-\lambda^{-2}\Delta_{\Gamma_k(x_n)}-c_k^{-2}\right)w_k,w_k\right\rangle_{L^2},$$
where $\Delta_{\Gamma_k(x_n)}$ denotes the (negative) Laplace-Beltrami operator on $\Gamma_k(x_n)$. Since $1-\eta_k$ is supported in the elliptic region $\{\zeta\in T^*\Gamma_k:\|\zeta\|>c_k^{-1}\}$, taking $\delta>0$ small enough we can arrange that on supp$\,\psi(1-\eta_k)$  the principal symbol of the operator $-\lambda^{-2}\Delta_{\Gamma_k(x_n)}-c_k^{-2}$ (considered as a semi-classical differential operator with a small parameter $\lambda^{-1}$) is lower bounded by a constant $C>0$ times the principal symbol of $-\lambda^{-2}\Delta_{\Gamma_k(x_n)}+1$. Therefore, by G\"arding's inequality we conclude
$$M\ge C\|w_k\|_{H^1(\Omega_k\setminus\Omega_{k-1})}^2,\quad C>0.\eqno{(4.27)}$$
On the other hand, by the trace theorem we have
$$\|g_k\|_{L^2(\Gamma_k)}^2\le C\lambda \|w_k\|_{H^1(\Omega_k\setminus\Omega_{k-1})}^2,\quad C>0.\eqno{(4.28)}$$
Clearly, (4.26) follows from (4.27) and (4.28).

Since $\Omega_k$ is strictly convex, the Neumann operator $N_k(\lambda)$ is a $\lambda-\Psi$DO with a principal symbol having a non-positive real part. The following properties of $N_k$ are proved in Section 3 of \cite{kn:CPV} (see Proposition 3.4).

\begin{lemma} There exists a constant $C>0$ such that we have
$${\rm Re}\,\left\langle N_k(\lambda)f,f\right\rangle_{L^2(\Gamma_k)}\le C\lambda^{-1/3}\|f\|_{L^2(\Gamma_k)}^2,\eqno{(4.29)}$$
$$\left\|[{\rm Op}_\lambda(\varphi),N_k(\lambda)]f\right\|_{L^2(\Gamma_k)}\le C\lambda^{-1/3}\|f\|_{H^1(\Gamma_k)}.\eqno{(4.30)}$$
\end{lemma}

Since $\|f_k\|_{H^1(\Gamma_k)}$ is equivalent to $\|g_k\|_{L^2(\Gamma_k)}$ and using the estimate
$$\|f_k\|_{H^1(\Gamma_k)}\le C\|f_k\|_{L^2(\Gamma_k)}+\|{\rm Op}_\lambda(1-\eta_k)f_k\|_{H^1(\Gamma_k)},$$
one can easily see that (4.24) follows from combining (4.25), (4.26), (4.29) and (4.30).
\eproof

Combining (4.23) and (4.24) and taking $\lambda$ big enough, we conclude
$$\|u\|_{H_k}+\|u_{k+1}\|_{L^2(K)}\le C\lambda^{-1}(\log\lambda)^{2^{k}}\left(\|v\|_{H_k}+\|v_{k+1}\|_{L^2({\bf R}^n\setminus\Omega_k)}\right).\eqno{(4.31)}$$
Clearly, (4.31) is equivalent to (3.2) for real $\lambda\gg 1$, which is the desired result.
\eproof

{\bf Acknowledgements.} A part of this work was carried out while F. Cardoso was visiting the University of Nantes in May 2009 with the support of the agreement Brazil-France in Mathematics - Proc. 69.0014/01-5. The first author is also partially supported by the CNPq-Brazil.

F. Cardoso

Universidade Federal de Pernambuco, 

Departamento de Matem\'atica, 

CEP. 50540-740 Recife-Pe, Brazil,

e-mail: fernando@dmat.ufpe.br\\

G. Vodev

Universit\'e de Nantes,

 D\'epartement de Math\'ematiques, UMR 6629 du CNRS,
 
 2, rue de la Houssini\`ere, BP 92208, 
 
 44332 Nantes Cedex 03, France,
 
 e-mail: vodev@math.univ-nantes.fr

\end{document}